\documentclass[10pt]{article}

\usepackage{amsmath, amssymb,amsfonts, amsopn, latexsym}

\newcommand{\U}{{\mathcal U}}
\newcommand{\0}{{\mathbf 0}}
\newcommand{\C}{{\mathbb C}}
\newcommand{\Z}{{\mathbb Z}}

\newcommand{\Q}{{\mathbb Q}}
\newcommand{\D}{{\mathbb D}} 
\newcommand{\N}{{\mathbb N}}
\newcommand{\cL}{{\mathbb L}}
\newcommand{\Proj}{{\mathbb P}}
\newcommand{\hyp}{{\mathbb H}}
\newcommand{\supp}{\operatorname{supp}}

\newcommand{\im}{\mathop{\rm im}\nolimits}
\newcommand{\rank}{\mathop{\rm rank}\nolimits}

\newcommand{\Adot}{\mathbf A^\bullet}
\newcommand{\Bdot}{\mathbf B^\bullet}

\newcommand{\strat}{{\mathfrak S}}

\newcommand{\cc}{{\operatorname{CC}}}
\newcommand{\ms}{{\operatorname{SS}}}

\newtheorem{defn0}{Definition}[section]
\newtheorem{prop0}[defn0]{Proposition}
\newtheorem{conj0}[defn0]{Conjecture}
\newtheorem{thm0}[defn0]{Theorem}
\newtheorem{lem0}[defn0]{Lemma}
\newtheorem{corollary0}[defn0]{Corollary}
\newtheorem{example0}[defn0]{Example}
\newtheorem{remark0}[defn0]{Remark}
\newtheorem{question0}[defn0]{Question}

\newenvironment{defn}{\begin{defn0}\hskip -.06in .}{\end{defn0}}

\newenvironment{thm}{\begin{thm0}\hskip -.06in .}{\end{thm0}}

\newenvironment{rem}{\begin{remark0}\hskip -.06in .\rm}{\end{remark0}}

\newcommand{\defref}[1]{Definition~\ref{#1}}

\newcommand{\thmref}[1]{Theorem~\ref{#1}}

\newcommand{\secref}[1]{Section~\ref{#1}}
\newcommand{\remref}[1]{Remark~\ref{#1}}

\newcommand{\qed}{\mbox{$\Box$}}
\newenvironment{proof}{\noindent {\bf Proof.}}{\qed\vskip 6pt}

\title{Vanishing Vanishing Cycles}

\author{David B. Massey}

\date{}

\begin{document}

\baselineskip = 14pt

\maketitle

\begin{abstract} If $\Adot$ is a bounded, constructible complex of sheaves on a complex analytic space $X$, and $f:X\rightarrow\C$ and $g:X\rightarrow\C$ are complex analytic functions, then the iterated vanishing cycles $\phi_g[-1](\phi_f[-1]\Adot)$ are important for a number of reasons.  We give a formula for the stalk cohomology $H^*(\phi_g[-1]\phi_f[-1]\Adot)_x$ in terms of relative polar curves, algebra, and the normal Morse data and micro-support of $\Adot$.  \end{abstract}

\sloppy




\section{Introduction}\label{sec:intro}

\footnotetext[1]{AMS subject classifications: 32B15, 32C35, 32C18, 32B10\newline Keywords: vanishing cycles, relative polar curve, characteristic cycle, discriminant}

Iterated vanishing and nearby cycles have appeared previously in the literature, notably in the work of M. Saito in \cite{mixedhodgemod}, Sabbah in \cite{sabbahprox}, and, in the motivic setting, in the work of Guibert, Loeser, and Merle in \cite{glm}. Nonetheless, such constructions and techniques may appear to be too abstract to be really useful in the down-to-Earth study of singular spaces. We wish to explain how iterated vanishing cycles arise
 in the study of the topology and geometry of singular spaces and functions on those spaces, even if one is interested solely in cohomology with $\Z$ coefficients.  

Suppose $X$ is a complex analytic space, and $Y$ is a closed analytic subspace of $X$. Let $j:Y\hookrightarrow X$ and $i:X-Y\hookrightarrow X$ denote the inclusions. Then, the hypercohomology modules $\hyp^k(X; \Z^\bullet_X)$, $\hyp^k(X; i_*i^*\Z^\bullet_X)$, $\hyp^k(X; j_*j^*\Z^\bullet_X)$, $\hyp^k(X; i_!i^!\Z^\bullet_X)$, and $\hyp^k(X; j_!j^!\Z^\bullet_X)$ are isomorphic to the following ordinary integral cohomology modules: $H^k(X; \Z)$, $H^k(X-Y; \Z)$, $H^k(Y, \Z)$, $H^k(X, Y; \Z)$, and $H^k(X, X-Y; \Z)$, respectively. In addition, if $f:X\rightarrow\C$ is a complex analytic function, then the complexes of nearby cycles, $\psi_f\Z^\bullet_X$, and vanishing cycles, $\phi_f\Z^\bullet_X$, are important for, at each $x\in V(f):=f^{-1}(0)$, the stalk cohomologies are, respectively, the integral cohomology and reduced integral cohomology of the Milnor fiber $F_{f, x}$ of $f$ at $x$. The language and results of the derived category allow us to treat all of these cohomology modules in a unified, elegant fashion. As general references, we recommend \cite{kashsch}, \cite{dimcasheaves}, \cite{schurbook}, and Appendix B of \cite{numcontrol}.

Characteristic cycles are fairly coarse data that one can associate to a complex of sheaves. Intuitively, they tell one the directions in which the Euler characteristic of the cohomology changes, and what that change is. If one wishes to retain more data from the complex, then, instead of considering the change in Euler characteristic, one can consider the cohomological {\it Morse modules} in each degree (see below). One should think of the Morse modules as telling one ``cohomological attaching data'' associated to how larger-dimensional strata are attached to lower-dimensional strata, with coefficients in the given complex of sheaves. The characteristic cycle is topologically important for it is related to the local Euler obstruction (see \cite{BDK}), the absolute polar varieties and local Morse inequalities (see Theorem 7.5 and Corollary 5.5 of \cite{numinvar}), and index formulas for the vanishing cycles (see \cite{ginsburg}, \cite{sabbahquel}, \cite{leconcept}, and below). In addition, the characteristic cycle of the intersection cohomology complex is of great importance in representation theory (see \cite{kazhdanlusztigtop} and \cite{bradenpams}).

The Morse modules, or coefficients of the characteristic cycle, of a complex $\Adot$ can be described in terms of vanishing cycles of $\Adot$ along functions with complex non-degenerate critical points. This means that the Morse modules, or coefficients of the characteristic cycle, of the complex of vanishing cycles $\phi_f\Adot$, contain iterated vanishing cycle data of the form $\phi_g\phi_f\Adot$, where $g$ is generic in some sense. Even for generic linear $g$, a formula for the Euler characteristic of the stalk cohomology is not well-known -- though it has appeared, in a slightly different form, in Corollary 4.6 of  \cite{hypercohom} (the last formula in the statement), and follows quickly from Theorems 3.3 and 5.5 of \cite{ginsburg} and Theorem 3.4.2 of \cite{bmm}. 

\medskip

The formula that we derive for generic linear $g$ led us to ask: what happens if $g$ is not so generic? The answer to this question is the main result, \thmref{thm:main}, of this short paper. Not surprisingly, the main result tells us that the computation of the iterated vanishing cycles boils down to calculating the relative polar curve of $f$ and $g$, and the intersection multiplicities of its components with $V(f)$ and $V(g)$. Even though the general flavor of the theorem is not surprising, the precise result is not at all obvious.

\bigskip

We will now present some technical background material, and state our result in a precise fashion.

\bigskip

Let $\U$ be an open neighborhood of the origin in $\C^{n+1}$, and let $\tilde f:(\U, \0)\rightarrow(\C, 0)$ be a complex analytic function. Suppose that $X$ is a complex analytic subspace of $\U$, and, for convenience, assume that $\0\in X$. Let $f$ denote the restrictions of $\tilde f$ to $X$. We write $V(f)$ for $f^{-1}(0)$.

Let $\strat$ be a complex analytic Whitney stratification of $X$ with connected strata. For $S\in\strat$, we let $\N_S$ and $\cL_S$ denote, respectively, the normal slice and complex link of the stratum $S$; see \cite{stratmorse}. We let $d_S:=\dim S$. Let $R$ be a regular, Noetherian ring with finite Krull dimension (e.g., $\Z$, $\Q$, or $\C$). Let $\Adot$ be a bounded complex of sheaves of $R$-modules on $X$, whose cohomology is constructible with respect to $\strat$.

For each $S\in\strat$ and $k\in\Z$, define the {\it degree $k$ Morse module of $S$ with respect to $\Adot$} to be the hypercohomology $m^k_S(\Adot):= \hyp^{k-d_S}(\mathbb N_S,\mathbb L_S; \Adot)=\hyp^{k}(\mathbb N_S,\mathbb L_S; \Adot[-d_s])$. If $S$ is the point-stratum $\{\0\}$, then the Morse module is isomorphic to the stalk cohomology of the vanishing cycles along a generic linear form; that is, if $\mathfrak l$ is a generic linear form, then $m^k_{\{\0\}}(\Adot)\cong H^k(\phi_{\mathfrak l}[-1]\Adot)_\0$. More generally, to describe the Morse modules of a higher-dimensional stratum $S$ in terms of vanishing cycles, one simply takes a normal slice to the stratum, to cut it down to an isolated point, restricts the complex $\Adot$ to the slice, shifts the restricted complex by $-d_S$, and then takes the (shifted by $-1$) vanishing cycles along a generic (affine) linear form on the slice.

We say that a stratum $S\in\strat$ is {\it $\Adot$-visible} if $m^*_S(\Adot)\neq 0$, and we let $\strat(\Adot)$ denote the set of $\Adot$-visible strata and let $\strat_\0(\Adot)$ denote the strata in $\strat(\Adot)$ which contain $\0$ in their closures. It is a theorem, Theorem 4.13 of \cite{micromorse}, that the {\it microsupport, $\ms(\Adot)$, of $\Adot$} (see \cite{kashsch}) is equal to $\bigcup_{S\in\strat(\Adot)}\overline{T^*_S\U}$.

If $R$ is an integral domain, define 
$$c_S(\Adot):=\sum_{k\in\Z} (-1)^{k}\rank(m^k_S(\Adot)),
$$
and define the {\it characteristic cycle of $\Adot$}, $\cc(\Adot)$, to be the analytic cycle in the cotangent space $T^*\U$ given by
$$
\cc(\Adot):=\sum_{S\in\strat}c_S(\Adot)\big[\overline{T^*_S\U}\big],
$$
where $T^*_S\U$ is the conormal space to $S$ in $\U$. (The reader should be aware that there are two or three different definitions of the characteristic cycle, differing by sign; see, for instance, \cite{kashsch} and \cite{schurbook}.  Previously, we have used a different convention, but also different notation for the characteristic cycle. The change in notation should help eliminate confusion for readers of our past work.) Throughout this paper, whenever we make a statement involving the characteristic cycle, we are assuming that $R$ is an integral domain, even though it will {\bf not} be explicitly mentioned.

The underlying set, $\left|\cc(\Adot)\right|$, is the {\it characteristic variety} of $\Adot$. The {\it characteristic cycle of $X$} is defined to be the characteristic cycle of $\Z^\bullet_X$.

\medskip

By taking a normal slice to a stratum, the calculation of the coefficient of $\overline{T^*_{S}\U}$ in $\cc(\Adot)$ is reduced to the calculation of the coefficient of the closure of the conormal space to a point-stratum. We use the origin as a convenient point and point out that, by definition, $c_\0(\Adot)$, the coefficient of $T^*_\0\U$ in $\cc(\Adot)$ is given in terms of the Euler characteristic of the stalk cohomology of the vanishing cycles by $c_\0(\Adot)=\chi(\phi_{\mathfrak l}[-1]\Adot)_\0$, where $\mathfrak l$ is a generic linear form.

\smallskip

Let $\im d\tilde f$ denote the image of $d\tilde f$, and let us recall the aforementioned index formulas for $\phi_f[-1]\Adot$, conjectured by Deligne, and proved independently by Ginsburg, Sabbah, and L\^e:

\begin{thm}\label{thm:lesabgin} {\rm (\cite{ginsburg}, \cite{sabbahquel}, \cite{leconcept})} Suppose that $\0$ is an isolated point in the support of $\phi_{f}[-1]\Adot$. 

Then, $(\0, d_\0\tilde f)$ is an isolated point in the intersection $|\cc(\Adot)|\cap \im d\tilde f$, and the Euler characteristic of the stalk cohomology of the vanishing cycles of $f$ is related to the intersection multiplicity of $\cc(\Adot)$ and image of $d\tilde f$ by $$\chi(\phi_{f}[-1]\Adot)_\0= \Big(\cc(\Adot)\cdot \operatorname{im}d\tilde f\Big)_{(\0, d_\0\tilde f)} = \sum_{S\in\strat_\0(\Adot)}c_S(\Adot)\,\left(\overline{T^*_{{}_{S}}X}\ \cdot\ \im d \tilde f\right)_{(\0, d_\0 \tilde f)}.$$
\end{thm}

\bigskip

Now, suppose that one wants the stalk cohomology modules of $\phi_{f}[-1]\Adot$, not merely the Euler characteristic. Then, one must begin with analogous data for $\Adot$; that is, one needs not simply $\cc(\Adot)$, but rather the closures of conormals spaces to strata together with the Morse modules of the strata with coefficients in $\Adot$.

\medskip

The refinement of \thmref{thm:lesabgin} that we proved in Theorem 5.3 of \cite{micromorse} was:

\begin{thm}\label{thm:isoimprov}{\rm (\cite{micromorse})} The origin is an isolated point in the support of $\phi_{f}[-1]\Adot$ if and only if $(\0, d_\0\tilde f)$ is an isolated point in the intersection $\ms(\Adot)\cap \im d\tilde f$, and if these equivalent conditions hold, then, for all $k$,
 $$H^{k}(\phi_f[-1]\Adot)_\0\ \cong \ \bigoplus_{S\in\strat_\0(\Adot)}\big(m^k_S(\Adot)\otimes_{{}_R}R^{{\delta_S}}\big),$$
where $\delta_S:=\left(\overline{T^*_{{}_{S}}X}\ \cdot\ \im d \tilde f\right)_{(\0, d_\0 \tilde f)}$.
\end{thm}

As we showed in Theorem 3.2 of \cite{hypercohom}, the intersection number $\delta_S$ of \thmref{thm:isoimprov} can be calculated in terms of the relative polar curve (see  \cite{hammlezariski}, \cite{teissiercargese}, \cite{leattach}, \cite{letopuse}, and \secref{sec:polar}), $\Gamma^1_{f,\mathfrak l}(S):=\Gamma^1_{f_{|_{\overline{S}}},\mathfrak l}$, where $\mathfrak l$ is a generic linear form, and $f_{|_S}$ is not constant. The formula that one obtains for a stratum $S$ of dimension at least $1$ is
$$
\delta_S= \big(\Gamma^1_{f,\mathfrak l}(S)\cdot V(f))_\0-  \big(\Gamma^1_{f,\mathfrak l}(S)\cdot V(\mathfrak l))_\0.
$$
In the case above, where $(\0, d_\0\tilde f)$ is an isolated point in the intersection $\ms(\Adot)\cap \im d\tilde f$, it is immediate that, if $S\in\strat_\0(\Adot)$ and $d_S\geq 1$, then $f_{|_S}$ is not constant. However, in the case below, it will be useful to let $\strat^f(\Adot)$ denote the set of strata in $\strat(\Adot)$ on which $f$ is not constant and let $\strat^f_\0(\Adot):=\strat^f(\Adot)\cap \strat_\0(\Adot)$, i.e., the set of $\Adot$-visible strata whose closures contain $\0$ but which are not contained in $V(f)$.

\bigskip

How can one generalize \thmref{thm:lesabgin} and \thmref{thm:isoimprov} to the case where the support of $\phi_f[-1]\Adot$ is of arbitrary dimension? Suppose that $\mathfrak l$ is a generic linear form (we shall not distinguish notationally between $\mathfrak l$ and ${\mathfrak l}_{|_{V(f)}}$). If $\0$ is an isolated point in the support of $\phi_f[-1]\Adot$, then it is trivial that there is an isomorphism of stalk cohomologies $H^k(\phi_{\mathfrak l}[-1]\phi_f[-1]\Adot)_\0\cong H^k(\phi_f[-1]\Adot)_\0$. 

Thus,  for generic linear $\mathfrak l$, \thmref{thm:lesabgin} and \thmref{thm:isoimprov} can be viewed as statements about $\chi(\phi_{\mathfrak l}[-1]\phi_f[-1]\Adot)_\0=c_\0(\phi_f[-1]\Adot)$ and $H^k(\phi_{\mathfrak l}[-1]\phi_f[-1]\Adot)_\0$, where $\delta_S$ is now {\bf defined} to be the value of $\big(\Gamma^1_{f,\mathfrak l}(S)\cdot V(f))_\0- \big(\Gamma^1_{f,\mathfrak l}(S)\cdot V(\mathfrak l))_\0$\, provided that $S\in\strat^f_\0(\Adot)$. 

\medskip

With these interpretations, we obtain a generalization of \thmref{thm:lesabgin} and \thmref{thm:isoimprov}; we prove as part of \thmref{thm:main} in this paper:

\medskip

\noindent{\bf Theorem: generic linear case}. {\it Let $\mathfrak l$ be a generic linear form and, for all $S\in\strat^f_\0(\Adot)$, let $\delta_S$  be the value of $\big(\Gamma^1_{f,\mathfrak l}(S)\cdot V(f))_\0- \big(\Gamma^1_{f,\mathfrak l}(S)\cdot V(\mathfrak l))_\0$.

Then, 
$$c_\0(\phi_f[-1]\Adot)=c_\0(\Adot)+\sum_{S\in\strat^f_\0(\Adot)}\delta_S \cdot c_S(\Adot).
$$

In fact, for all $k\in\Z$,
$$
H^k(\phi_{\mathfrak l}[-1]\phi_f[-1]\Adot)_\0\cong m^k_{\{\0\}}(\Adot)\ \oplus\bigoplus_{S\in\strat^f_\0(\Adot)}\big(m^k_S(\Adot)\otimes R^{\delta_S}\big).
$$
}

As we discuss in \remref{rem:braden}, the result above allows us to quickly obtain Theorem 1 of \cite{bradenpams}.

\bigskip

Seeing that the generic linear form case of our theorem provides a formula for the stalk cohomology of the iterated vanishing cycles $\phi_{\mathfrak l}[-1]\phi_f[-1]\Adot$ in terms of Morse modules and intersection multiplicities with relative polar curves, one is naturally led to ask: is there a more general formula where all occurrences of $\mathfrak l$ are replaced by a more general function $g$? The answer is: yes.
\smallskip

Suppose that we give ourselves a second complex analytic function $\tilde g:(\U, \0)\rightarrow(\C, 0)$, and let $g$ denote the restriction of $\tilde g$ to $X$. For each stratum $S\in\strat^f(\Adot)$, we define in \secref{sec:polar}, following our work in \cite{enrichpolar}, a {\it relative polar set} $|\Gamma_{f, \tilde g}(S)|$ and, if this set is $1$-dimensional, we define a corresponding {\it relative polar curve} (as a cycle) $\Gamma^1_{f, \tilde g}(S)$; this relative polar curve agrees with the traditional one in the case where $\tilde g$ is a generic linear form.

We also define $\widehat\Gamma^1_{f, \tilde g}(S)$ to be the sum of the components of the cycle $\Gamma^1_{f, \tilde g}(S)$ which contain the origin, but which are not contained in $V(g)$. We define $\big|\Gamma_{f,\tilde g}(\Adot)\big|:=\bigcup_{S\in\strat^f(\Adot)}\big|\Gamma_{f,\tilde g}(S)\big|$. For all $S\in\strat(\Adot)$, we define
  $$\hat\delta_{f, \tilde g}(S):=\sum_{C\text{ comp. of } \widehat\Gamma^1_{f,\tilde g}(S)}\Big( \big(C\cdot V(f)\big)_\0-\operatorname{min}\big\{ \big(C\cdot V(f)\big)_\0, \big(C\cdot V(g)\big)_\0\big\}\Big).
 $$
In the case where $\tilde g$ is a generic linear form, $\hat\delta_{f, \tilde g}(S)$ will agree with $\delta_S$ from generic linear form case of our theorem.

\bigskip

Our general result is \thmref{thm:main}:

\smallskip

\noindent {\bf Theorem}. {\it Suppose that $\dim_\0 V(f)\cap \big|\Gamma_{f,\tilde g}(\Adot)\big|\leq 0$.

Then, for generic $[\sigma : \lambda]\in\Proj^1$, for all $k\in\Z$,
$$
H^k(\phi_g[-1]\phi_f[-1]\Adot)_\0\cong H^k(\phi_{\sigma f+\lambda g}[-1]\Adot)_\0\ \oplus\bigoplus_{S\in\strat^f(\Adot)}\big(m^k_S(\Adot)\otimes R^{\hat\delta_{f, \tilde g}(S)}\big).
$$
}

If $\tilde g$ is a generic linear form $\mathfrak l$, then $H^k(\phi_{\sigma f+\lambda g}[-1]\Adot)_\0\cong H^k(\phi_{\mathfrak l}[-1]\Adot)_\0\cong m^k_{\{\0\}}(\Adot)$ and $\hat\delta_{f, \tilde g}(S)$ is equal to $\delta_S$; thus, the main theorem generalizes the result in the generic linear form case.

Pencils of Milnor fibrations of the form $\sigma f+\lambda g$ have been studied in detail by Caubel in \cite{caubelthesis} and \cite{caubel}. Caubel's main technique, the ``tilting in the Cerf diagram'', first used by L\^e and Perron in \cite{leperron}, was the essential idea that we used to {\bf see} that the main theorem of this paper is true, but the technical details needed to make the tilting argument rigorous are formidable; the proof that we give here uses our technical result on the enriched discriminant from \cite{enrichpolar}.

\section{The General Relative Polar Curve}\label{sec:polar}

We continue with all of the notation from \secref{sec:intro}. Our treatment of the relative polar curve follows \cite{enrichpolar}, except that we avoid discussing enriched cycles.

\bigskip

Suppose that $M$ is a complex submanifold of $\U$. Recall:

\begin{defn}\label{def:relconorm} The relative conormal space $T^*_{\tilde f_{|_M}}\U$ is given by 
$$
T^*_{\tilde f_{|_M}}\U :=\{(x, \eta)\in T^*\U\ |\ \eta(T_xM\cap \ker d_x\tilde f)=0\}.
$$

If $M\subseteq X$, then $T^*_{\tilde f_{|_M}}\U$ depends on $f$, but not on the particular extension $\tilde f$. In this case, we write $T^*_{f_{|_M}}\U$ in place of $T^*_{\tilde f_{|_M}}\U$.
\end{defn}

\smallskip

Let $\pi:T^*\U\rightarrow\U$ denote the projection. 

\bigskip

We are going to define a relative polar set $\big|\Gamma_{f,\tilde g}(S)\big|$ and a relative polar cycle  $\Gamma^1_{f,\tilde g}(S)$, as we did in \cite{enrichpolar}. If $\tilde g$ is a generic linear form, and $f_{|_S}$ is not constant, it is easy to show that our $\big|\Gamma_{f,\tilde g}(S)\big|$ is purely $1$-dimensional and that $\Gamma^1_{f,\tilde g}(S)$ is reduced, and agrees with all of the definitions/characterizations of the relative polar curve used in \cite{hammlezariski}, \cite{teissiercargese}, \cite{leattach}, \cite{letopuse} by Hamm, L\^e, and Teissier. The point of \defref{def:polarcurve} is that it seems to be the ``correct'' definition of the relative polar curve even when $\tilde g$ is not so generic. Note that, in the traditional case where $\tilde g$ is a non-zero linear form, $d_x\tilde g$ is a ``constant'' non-zero covector.

\begin{defn}\label{def:polarcurve} If $S\in\strat$ and $f_{|_S}$ is not constant, we define the  {\bf relative polar set}, $\big|\Gamma_{f,\tilde g}(S)\big|$, to be $\pi\left(\overline{T^*_{f_{|_S}}\U}\ \cap\ \im d\tilde g\right)$. The  {\bf relative polar set}, $\big|\Gamma_{f,\tilde g}(\Adot)\big|$, is defined by
$$\big|\Gamma_{f,\tilde g}(\Adot)\big|:=\bigcup_{S\in\strat^f(\Adot)}\big|\Gamma_{f,\tilde g}(S)\big|.$$

If $C$ is a (reduced) $1$-dimensional component of $\big|\Gamma_{f,\tilde g}(S)\big|$, then $\overline{T^*_{f_{|_S}}\U}$ and $\im d\tilde g$ intersect properly along a unique $1$-dimensional component $C^\prime:=d\tilde g(C)=\{(x, d_x\tilde g)\ |\ x\in C\}$ such that $\pi(C^\prime)=C$. Thus, the intersection number $p_C(S):= \left(\left[\overline{T^*_{f_{|_S}}\U}\right]\cdot[\im d\tilde g]\right)_{C^\prime}$ is well-defined.

In particular, if $\big|\Gamma_{f,\tilde g}(S)\big|$ is purely $1$-dimensional, then we may, and do, define the {\bf relative polar curve}, $\Gamma^1_{f,\tilde g}(S)$, to be the properly pushed-forward cycle 
$$\pi_*\left(\left[\overline{T^*_{f_{|_S}}\U}\right]\cdot [\im d\tilde g]\right)=\sum_Cp_C(S)[C],$$
where the sum is over all of the components of $\big|\Gamma_{f,\tilde g}(S)\big|$.

If $\big|\Gamma_{f,\tilde g}(\Adot)\big|$ is purely $1$-dimensional (respectively, is $1$-dimensional at the origin), then we define the {\bf relative polar curve, as a cycle (respectively, as a cycle germ at the origin)} to be 
$$
\Gamma^1_{f,\tilde g}(\Adot):=\sum_{S\in\strat^f(\Adot)}\Gamma^1_{f,\tilde g}(S).
$$

\end{defn}

\begin{rem} It is trivial to show that $\big|\Gamma_{f,\tilde g}(S)\big|$ has no zero-dimensional components. Thus, using the convention that the empty set has dimension $-\infty$, the condition that $\dim_\0\big|\Gamma_{f,\tilde g}(\Adot)\big|\leq 1$ is equivalent to $\big|\Gamma_{f,\tilde g}(\Adot)\big|$ being purely one-dimensional at $\0$.

It will also be important to us later what happens in the case of a $1$-dimensional stratum $S$ along which $f$ is not constant. It is trivial to see that, in this case, the cycle $\Gamma^1_{f,\tilde g}(S)$ is reduced and simply equals  $[\overline{S}]$.
\end{rem}

\smallskip

\begin{defn}\label{def:delta} If $\dim_\0 V(f)\cap \big|\Gamma_{f,\tilde g}(\Adot)\big|\leq 0$ and $\dim_\0 V(g)\cap \big|\Gamma_{f,\tilde g}(\Adot)\big|\leq 0$, then, for all $S\in\strat^f(\Adot)$, we define $\delta_{f, \tilde g}(S)$ to be the difference of intersection numbers $\big(\Gamma^1_{f,\tilde g}(S)\cdot V(f)\big)_\0-\big(\Gamma^1_{f,\tilde g}(S)\cdot V(g)\big)_\0$.
\end{defn}

\begin{rem} Note that the two dimension hypotheses of \defref{def:delta} are satisfied when $\tilde g$ is a generic linear form (see \cite{enrichpolar}, Proposition 3.13). Also, by the work of Hamm, L\^e, and Teissier, if $\mathfrak l$ is a generic linear form, then $\Gamma^1_{f,\mathfrak l}(S)$ is reduced and $\big(\Gamma^1_{f, \mathfrak l}(S)\cdot V(\mathfrak l)\big)_\0$  is the multiplicity of $\Gamma^1_{f,\mathfrak l}(S)$ at the origin; thus, $\delta_{f, \tilde g}(S)\geq 0$, when $\tilde g$ is a generic linear form.
\end{rem}

\section{The Vanishing Vanishing Theorem}\label{sec:main}

We continue using all of the notation from the previous two sections.

\smallskip

 Suppose that $\dim_\0\big|\Gamma_{f,\tilde g}(\Adot)\big|\leq 1$. If $S\in\strat^f_\0(\Adot)$, we let $\big|\widehat\Gamma_{f,\tilde g}(S)\big|$ denote the union of the components of $\big|\Gamma_{f,\tilde g}(S)\big|$ which contain the origin, but are not contained in $V(g)$; let $\widehat\Gamma^1_{f,\tilde g}(S)$ denote the corresponding cycle. Let $\big|\widehat\Gamma_{f,\tilde g}(\Adot)\big|:=\bigcup_{S\in\strat^f_\0(\Adot)}\big|\widehat\Gamma_{f,\tilde g}(S)\big|$. Note that Lemma 3.10 of \cite{enrichpolar} implies that no component of $\big|\widehat\Gamma_{f,\tilde g}(\Adot)\big|$ is contained in $V(f)$.
 
 Let $C$ denote a possibly non-reduced component of the cycle $\widehat\Gamma^1_{f,\tilde g}(\Adot)$; we will write $|C|$ for the underlying analytic set of $C$ (i.e., $C$ with its reduced structure). Thus, 
 $$
 C=\Big(\sum_{S\in\strat^f(\Adot)}p_{|C|}(S)\Big)|C|,
 $$
 where $p_{|C|}(S)$ is as in \defref{def:polarcurve}.
 
 For all $S\in\strat(\Adot)$, define
  $$\hat\delta_{f, \tilde g}(S):=\sum_{C\text{ comp. of } \widehat\Gamma^1_{f,\tilde g}(S)}\Big( \big(C\cdot V(f)\big)_\0-\operatorname{min}\big\{ \big(C\cdot V(f)\big)_\0, \big(C\cdot V(g)\big)_\0\big\}\Big).
 $$
 Note that, in the formula above, if $\0$ is not in the closure of $S$, then $\hat\delta_{f, \tilde g}(S)$ is automatically equal to $0$.
 
\bigskip

We will now prove the new result of this paper. Our method is to analyze the situation using the discriminant and Cerf diagram, the same method used first by L\^e in \cite{leattach}, and subsequently in many works by various authors, such as L\^e and Perron in \cite{leperron}, Tib\u ar in \cite{tibarthesis}, \cite{tibarcar}, and \cite{tibarvan}, and  Caubel in his thesis \cite{caubelthesis} and \cite{caubel}.  The main technical difficulty is that one needs a derived category version of the discriminant and/or Cerf diagram. However, we proved in \cite{enrichpolar} that such a discriminant exists.

\begin{thm}\label{thm:main} Suppose that $\dim_\0 V(f)\cap \big|\Gamma_{f,\tilde g}(\Adot)\big|\leq 0$.

Then, for generic $[\sigma : \lambda]\in\Proj^1$, for all $k\in\Z$,
$$
H^k(\phi_g[-1]\phi_f[-1]\Adot)_\0\cong H^k(\phi_{\sigma f+\lambda g}[-1]\Adot)_\0\ \oplus\bigoplus_{S\in\strat^f_\0(\Adot)}\big(m^k_S(\Adot)\otimes R^{\hat\delta_{f, \tilde g}(S)}\big).
$$

\medskip

In particular, suppose that $\mathfrak l$ is a linear form, which is generic with respect to the fixed $f$ and $\strat$. For all $S\in\strat^f_\0(\Adot)$, let $\delta_S:=\big(\Gamma^1_{f,\mathfrak l}(S)\cdot V(f))_\0- \big(\Gamma^1_{f,\mathfrak l}(S)\cdot V(\mathfrak l))_\0$.

Then, 
$$c_\0(\phi_f[-1]\Adot)=c_\0(\Adot)+\sum_{S\in\strat^f_\0(\Adot)}\delta_S \cdot c_S(\Adot).
$$

In fact, for all $k\in\Z$,
$$
H^k(\phi_{\mathfrak l}[-1]\phi_f[-1]\Adot)_\0\cong m^k_{\{\0\}}(\Adot)\ \oplus\bigoplus_{S\in\strat^f_\0(\Adot)}\big(m^k_S(\Adot)\otimes R^{\delta_S}\big).
$$

\end{thm}
\begin{proof} We first prove the generic linear form case, and then use that to prove the general case.

\smallskip

\noindent{Generic linear form case}: 

The statement about the coefficients in the characteristic cycle is proved in Corollary 4.6 of  \cite{hypercohom} (the last formula in the statement); however, it also follows quickly from Theorems 3.3 and 5.5 of \cite{ginsburg} and Theorem 3.4.2 of \cite{bmm}.

If the base ring $R$ is a field, the stalk cohomology statement follows at once by combining the characteristic cycle coefficient formula with Lemma 2.3 of \cite{micromorse}, where we used perverse cohomology to extract individual Betti numbers from characteristic cycle formulas. That is, in the first formula of the above theorem, we replace $\Adot$ by the perverse cohomology ${}^{\mu}\negmedspace H^k(\Adot)$ (see \cite{bbd}, \cite{kashsch}, or, for a quick summary, \cite{micromorse}), and use that $$\cc({}^{\mu}\negmedspace H^k(\Adot))=\sum_{S\in\strat} b_{{}_{k-d_S}}(\N_S, \cL_S;\Adot)\left[\overline{T^*_S\U}\right],$$ 
where $b_{{}_{j}}(\N_S, \cL_S;\Adot)$ denotes the $j$-th relative Betti number, i.e., $b_{{}_{j}}(\N_S, \cL_S;\Adot):=\dim_{{}_R}\hyp^j(\N_S, \cL_S;\Adot)$.

For general $R$, the isomorphism is obtained by ``enriching'' the proof of Lemma 3.7 in \cite{numcontrol}. That is, one uses precisely the proof of Lemma 3.7 in \cite{numcontrol}, except that one replaces cycles by enriched cycles and ordinary intersection theory by enriched intersection theory; such enriched proofs are discussed in \cite{singenrich}.

\bigskip

\noindent{General case}:

The idea of this proof is simple: push the complex $\Adot$ down onto an open polydisk in $\C^2$ via the map $(g, f)$, and then apply the generic linear case to the pushed-forward complex. The details are far from simple. However, the major technical piece of the proof -- that there is a well-behaved derived category version of the discriminant -- appears in our work \cite{enrichpolar}.

\medskip

For $\epsilon, \delta, \rho>0$, let $N^\epsilon_{\delta, \rho}:=B_\epsilon\cap g^{-1}({\stackrel{\circ}{\D}}_{{}_\delta})\cap f^{-1}({\stackrel{\circ}{\D}}_{{}_\rho})$. Let $(\Adot)^\epsilon_{\delta, \rho}$ be the restriction of $\Adot$ to $N^\epsilon_{\delta, \rho}$, and let $T^\epsilon_{\delta, \rho}$ be the restriction of the map $(g, f)$ to a (proper) map from $N^\epsilon_{\delta, \rho}$ to ${\stackrel{\circ}{\D}}_{{}_\delta}\times {\stackrel{\circ}{\D}}_{{}_\rho}$.

As we proved in the Main Theorem of \cite{enrichpolar}, for all sufficiently small $\epsilon>0$, there exist $\delta, \rho>0$ such that the derived push-forward $\Bdot:=R(T^\epsilon_{\delta, \rho})_*(\Adot)^\epsilon_{\delta, \rho}$ is complex analytically constructible with respect to the stratification given by 
$$
\{{\stackrel{\circ}{\D}}_{{}_\delta}\times {\stackrel{\circ}{\D}}_{{}_\rho}-\Delta_{\Adot}(\tilde g, \tilde f),\ ({\stackrel{\circ}{\D}}_{{}_\delta}\times {\stackrel{\circ}{\D}}_{{}_\rho})\cap\Delta_{\Adot}(\tilde g, \tilde f)-\{\0\}, \ \{\0\}\},
$$
where $\Delta_{\Adot}(\tilde g, \tilde f)$ is the {\it $\Adot$-discriminant}, and is equal to 
$T^\epsilon_{\delta, \rho}\Big(N^\epsilon_{\delta, \rho}\cap \left(\big|\Gamma_{f,\tilde g}(\Adot)\big|\cup\supp\phi_f[-1]\Adot\right)\Big)$. 

Using $(v, u)$ for coordinates on ${\stackrel{\circ}{\D}}_{{}_\delta}\times {\stackrel{\circ}{\D}}_{{}_\rho}$, it is important to realize that the visible components of $\Delta_{\Adot}(\tilde g, \tilde f)$, other than (possibly) ${\stackrel{\circ}{\D}}_{{}_\delta}\times\{0\}$,  form the reduced cycle $\Gamma_{u,v}(\Bdot)$. The  description of the Morse modules of $\Bdot$ is contained in the last paragraph of Corollary 4.4 of \cite{enrichpolar}; we describe a portion of this result here.

Let $P$ be an irreducible component of $\Delta_{\Adot}(\tilde g, \tilde f)$, other than $V(u)={\stackrel{\circ}{\D}}_{{}_\delta}\times\{0\}$. Then, corresponding to our previous $\delta_S$, for generic $[\sigma : \lambda]\in\Proj^1$, we let 
$$\delta_P \ = \  \big(P\cdot V(u))_\0-  \big(P\cdot V(\sigma u+\lambda v))_\0 \ = \ \big(P\cdot V(u)\big)_\0-\operatorname{min}\big\{ \big(P\cdot V(u)\big)_\0, \big(P\cdot V(v)\big)_\0\big\},
$$
where, if $P=V(v)=\{0\}\times{\stackrel{\circ}{\D}}_{{}_\rho}$, we mean that $\delta_P=0$.
Then, it follows immediately from Corollary 4.4 of \cite{enrichpolar}, together with the intersection number formula for proper push-forwards, that
$$
\bigoplus_{P\in\strat^u_\0(\Bdot)}\big(m^k_P(\Bdot)\otimes R^{\delta_P}\big) \ = \ \bigoplus_{S\in\strat^f_\0(\Adot)}\big(m^k_S(\Adot)\otimes R^{\hat\delta_{f, \tilde g}(S)}\big).
$$

Now, proper push-forwards commute with taking vanishing cycles (see \cite{kashsch}). Thus, we have
$$
\phi_{\sigma u+\lambda v}[-1]\Bdot \ \cong R(T^\epsilon_{\delta, \rho})_*\left(\phi_{\sigma f+\lambda g}[-1]\Adot\right).
$$
and
$$
\phi_v[-1]\phi_u[-1]\Bdot \ \cong R(T^\epsilon_{\delta, \rho})_*\left(\phi_g[-1]\phi_f[-1]\Adot\right).
$$

By taking $\epsilon$ small enough, constructibility tells us that 
$$
m^k_{\{\0\}}(\Bdot) \ \cong \ H^k(\phi_{\sigma u+\lambda v}[-1]\Bdot)_\0.
 \ \cong \ \hyp^k\left(\big(T^\epsilon_{\delta, \rho}\big)^{-1}(\0);\ \phi_{\sigma f+\lambda g}[-1]\Adot\right)  \ \cong
 $$
$$\hyp^k\left(B_\epsilon\cap V(\sigma f+\lambda g);\ \phi_{\sigma f+\lambda g}[-1]\Adot\right) \ \cong \ 
H^k(\phi_{\sigma f+\lambda g}[-1]\Adot)_\0.
$$
and
$$
H^k(\phi_{\sigma u+\lambda v}[-1]\phi_u[-1]\Bdot)_\0 \ \cong \ H^k(\phi_v[-1]\phi_u[-1]\Bdot)_\0
 \ \cong \ \hyp^k\left(\big(T^\epsilon_{\delta, \rho}\big)^{-1}(\0);\ \phi_g[-1]\phi_f[-1]\Adot\right)  \ \cong
 $$
$$\hyp^k\left(B_\epsilon\cap V(f, g);\ \phi_g[-1]\phi_f[-1]\Adot\right) \ \cong \ 
H^k(\phi_g[-1]\phi_f[-1]\Adot)_\0.
$$

The general result for $f$, $g$, and $\Adot$ now follows from the generic linear result for $u$ and $\Bdot$. \end{proof}

\bigskip

\begin{rem}  It is an important part of \thmref{thm:main} that the amount of genericity that we need for $g$ is precisely that $\dim_\0 V(f)\cap \big|\Gamma_{f,\tilde g}(\Adot)\big|\leq 0$. This dimension condition can, in fact, be checked in practice. In particular, we do {\bf not} need for the origin to be an isolated point in the support of $\phi_g[-1]\phi_f[-1]\Adot$ or of $\phi_{\sigma f+\lambda g}[-1]\Adot$. However, if the origin is, in fact, an isolated point in the support of $\phi_{\sigma f+\lambda g}[-1]\Adot$, then the stalk cohomology of $\phi_{\sigma f+\lambda g}[-1]\Adot$ at the origin can be calculated via \thmref{thm:isoimprov}.

We should remark that, as we showed in Theorem 3.14 of \cite{enrichpolar}, the condition that $\dim_\0 V(f)\cap \big|\Gamma_{f,\tilde g}(\Adot)\big|\leq 0$ is equivalent to the origin being an isolated point in the support of $\phi_g[-1]\psi_f[-1]\Adot$ (note the nearby cycles along $f$, not the vanishing cycles), or not in the support at all.
\end{rem}

\begin{rem}\label{rem:braden} Suppose that $\tilde f$ is itself a linear form, which is not necessarily generic, and that, at the origin, for generic linear $\mathfrak l$,  $\big|\Gamma_{f,\mathfrak l}(\Adot)\big|$ is a collection of lines; this would, for instance, be the case if $\Adot$ were constructible with respect to a conic stratification of affine space. Then, for generic linear $\mathfrak l$, for all $S\in \strat^f_\0(\Adot)$, $\delta_S=0$, and so 
$H^k(\phi_{\mathfrak l}[-1]\phi_f[-1]\Adot)_\0\cong H^k(\phi_{\mathfrak l}[-1]\Adot)_\0$ and $c_\0(\phi_f[-1]\Adot)= c_\0(\Adot)$. Thus, we recover Theorem 1 of \cite{bradenpams}.
\end{rem}

\newpage
\bibliographystyle{plain}
\bibliography{Masseybib}
\end{document}